\documentclass[12pt,fleqn]{article}
\begin{document}
\title{\textbf{Oscillating Population Models}}
\author{B S Lakshmi\\Center For Applicable Mathematics and computer Sciences\\B M Birla Science
Center\\Adarsh Nagar\\Hyderabad 500063}
\date{}
\maketitle

\begin{abstract}
Oscillating population model realistic situations in different
contexts.We examine this situation with reasonable mathematical
models and come to interesting conclusions,such as for
example,that the population at most points of the cycle
approximately equals half the maximum attainable
population.\end{abstract}
\section{Introduction}
 The Logistic equation
\begin{equation}\frac{d(P(t))}{dt}=r(M-P(t))P(t) \ , \
r>0\end{equation} where P(t) is the population at given time t and
M is the maximum sustainable population was formulated to refine
the well-known Malthusian Model of population growth.\cite{r1,r2}
This model takes into consideration that the resources available
are limited. The logistic form is a very useful mathematical tool
in describing the phenomena of growth.Experiments have shown a
fairly satisfactory fit between empirical observations and
theoretical calculations based on this logistic from of
Verlhurst.This logistic model gives a theory of growth in the
simplest stage --- a struggle for existence. The Logistic equation
can be treated as an equation of the Bernoulli type.By using a
partial fraction decomposition it becomes
\begin{equation}\frac{1}{M}(\frac{1}{P}+\frac{1}{M-P})\frac{dP}{dt}
= r \end{equation}On integration this yields \begin{equation}
\frac{P}{M-P} = Ae^{rMT}\end{equation} A, the constant of
integration is chosen from the initial condition\\ $P=P_{0}$ at \
$t=t_{0}$\\ Thus from (3) we get
\begin{equation}P(t)=\frac{MP_{0}}{P_{0}+(M-P_{0})e^{-rM(t-t_{0})}}
\end{equation}
We now consider equation(1) in the following context:\\ The
maximum attainable population would be different at different
periods.For example this would be the case when\\ seasonal factors
come into play as in the case of insect populations like
mosquitoes which manifest highs and lows at different times of the
year due to temperature dependent factors.\\Secondly the
well-known 10 year cycle in lynx and snow-shoe hare populations of
the boreal forests of Canada is perhaps controlled by the
interaction which can be best represented as
follows\cite{r3}:\\vegetation$\longrightarrow$ eaten by hares
$\longrightarrow$ eaten by lynx\\Thus the lynx exhibits ``driven
oscillations'' in its population.\\It maybe mentioned that
oscillating populations have also been studied in certain mammals
in a predator-prey model \cite{r4,r5}\section{Oscillating
Populations} We now model the above considerations in two simple
ways.\subsection{} First we consider the case where M assumes two different values\\
$M_{1}$ and $ M_{2}$,
\[
M=M_{1}\hspace{20mm} 0\leq t \leq h/2\]

\[
M=M_{2}\hspace{20mm} h/2\leq t \leq h
\]
We now write equation (1) as
\[
\frac{1}{r}\frac{dP}{dt} = MP-P^{2}=-(P-M/2)^{2}+\frac{M^{2}}{4}\]
 Putting \[ P-M/2=W\],we get the Riccati equation,
\begin{equation}\frac{1}{r}\frac{dW}{dT} = -W^{2}+ M^{2}/4-\frac{1}{2}\frac{dM}{dt}
\end{equation}
We observe that equation (5) shows  that P has a periodic solution
with the same period as M.This would also be expected from (1)
itself.

Let$ M_{1}$ and $M_{2}$ be the maximum sustainable populations
during different time periods.Then as indicated above equation (4)
would have two different forms:
\begin{equation}\
P_{1}(t) =
\frac{M_{1}P_{0}}{P_{0}+(M_{1}-P_{0})e^{-rM_{1}h/2}}\hspace{25mm}
0\leq t \leq h/2\end{equation} \begin{equation}
 P_{2}(t) =\frac{M_{2}P_{1}}{P_{1}+(M_{2}-P_{1})e^{-rM_{2}h/2}}\hspace{25mm}
h/2 \leq t \leq h\end{equation}Equations (6) and (7)share the
interlinkage of populations in the two intervals ---  this
ofcourse is quite expected.\\ For large positive values of $M_{1}$
and $M_{2}$ ,$M_{1}\gg P_{0}$ ,$M_{2}\gg P_{1}$ in equations
(6)and (7),we further get,
\[P_{1}(t)\approx M_{1 } \ \ \ \ \ P_{2}\approx M_{2}\] .This shows that the maximal values are soon
attained in the two intervals. If we now impose the condition
\[\frac{M_{1}+M_{2}}{2} \approx \langle P\rangle ,\] that is, that
on the average,the population in the interval $0<t<h$ is the mean
of the two maximum sustainable populations in the two
sub-intervals,then we deduce that
\[P_{1}\approx P_{2}\]
\subsection{}

 We now consider the case when M is periodic,with period h which models oscillating sustainability.
As already remarked after (5) or as can be seen from (1)\\ there
exists a solution P which is periodic with the same period as M.In
other words the population  could follow the same periodic pattern
as M.Let us first scale equation (1) by dividing both sides with a
large number $N^{2}$ so the maximum attainable population is
normalized to 1.\\Thus we have
\begin{equation}\frac{1}{N^{2}}\frac{dP}{dt} = \frac{M
P}{N^{2}}-\frac{P^{2}}{N^{2}}\end{equation}  Put \[ \frac{M}{N}
\equiv {\mathcal{M}} \  \ \frac{P}{N} \equiv {\mathcal{P}},\] and
integrating both sides from 0 to h, we get,
\[\int_{0}^{h}\frac{1}{N}d{\mathcal{P}} =
\int_{0}^{h}({\mathcal{MP}} - {\mathcal{P}}^{2})dt\] The left hand
side = 0, since ${\mathcal{P}}$ is taken to be periodic. The right
hand side can be written as
\[-\int_{0}^{h}[({\frac{\mathcal{M}}{2}-P})^{2}-\frac{{\mathcal{M}}^{2}}{4}
] dt\]whence we get
\[\int_{0}^{h}({\mathcal{M}-P}/2)^{2} \ dt = \int_{0}^{h}{\mathcal{M}}^{2}/4 \
dt\]As can be seen, for small ${\mathcal{M} ,{M}}^{2}$ can be
neglected except in a small interval near its maximum, whence
\[\int({\mathcal{M}}/2-{\mathcal{P}})^{2} \ dt \simeq 0,\]
\[ {\mathcal{M}}/2 \simeq {\mathcal{P}}\] This shows that the maximum
sustainable population  which continuously and periodically
changes with time is twice the actual population at that time at
any point in the cycle. \vspace{10mm}\\ \subsection{} \vspace{5mm}
We explicitly solve the  equation (1) to get a graphical feel of
the above considerations  .Putting \[z = \frac{1}{P} \ \ \ \ \ z'=
\frac{-P'}{P^{2}}\] ,we get
\begin{equation}-z' - rMz + r = 0 \end{equation} The integrating Factor
for (9) is $e^{r\int_{t_{0}}^{t}M(t)dt} \equiv Q $ \\Solving we
get
\[z = \frac{r}{Q}\int_{T}^{t}Q dt\] and hence \begin{equation} P =\frac{1}{z} =
\frac{Q}{r\int_{T}^{t}Q dt}\end{equation}
\[Q = e^{r\cos t_{0}}e^{-r\cos t}\]
(10) gives the population at time t.\cite{r1}.We now graphically
consider (10) for a simple special case:
\[M(t)= \sin t\]
\begin{equation}y = \frac{r}{e^{-r\cos t}}\int_{T}^{t}e^{-r\cos
t}dt \end{equation}

\section{Remarks} We remark that equation(1)
maybe considered with multiple parameters \cite{r6,r7}.On the
other hand,the behaviour of the solution is dependent on the
parameters.\cite{r8} the value of P typically increases and is
proportional to its initial value,but at large times it
$\rightarrow M$, here treated as a constant,this being independent
of the initial value.On the other hand the discrete version of the
logistic equation,equation(1) is even more complicated in that as
small r increases, first one, then two and subsequently more and
more number of solutions appear.P oscillates betweem them and
ultimately becomes totally random.
\section*{Acknowledgement}I am thankful to Dr.B.G.Sidharth for
useful discussions.

\end{document}